\documentclass[preprint,aps,showpacs]{revtex4}
\begin{document} 
%\draft 
\preprint{UQ Theory: December 1992}

\title{Approximate Dynamical Systems}
\author{K. R. W. Jones}
\address{Physics Department, University of Queensland,\\
St Lucia 4072, Brisbane, Australia.}
\date{December 1992}
\begin{abstract}
Working notes on setting up approximate dynamical systems
and nonlinear eigenvalue problems, here embedded within
the theory of complex nonlinear dynamics. Computations
parallel those of linear quantum theory except that we
use functional methods rather than Hilbert space.
\end{abstract} 
\bigskip        
%\pacs{}
\maketitle

\section{Introduction}
This covers basic concepts behind the theory of
{\em approximate dynamical systems}\cite{jones}, the 
theory of generalized quantum dynamics, functional 
methods, theory of propagators, and the basics of 
nonlinear spectral theory and nonlinear
functional analysis.

It begins with the key idea of the subject, namely the
use of {\em action principles} to define, and obtain
{\em approximate dynamical systems}. Subsequently
we see some examples of useful systems of this 
kind that generate nonlinear equations.

\section{Action principles and approximation theory}
What we really need, to make the business of practical
computations efficient, is a {\em larger mathematics}
which contains both exact and approximate equations
of motion --- treated in a unified manner. This is
best done working from the principle of least action.

At the deepest level, physical theories when considered
as dynamical systems, derive from action principles. It
easy to state these in complete generality, as a shell
into which we plug a Lagrangian and turn the crank. 
Here we concentrate upon {\em equations of motion} 
as ``solutions'' of the variational problem
\begin{equation}
\label{leastaction}
\frac{\delta}{\delta x(t)}
\int_{t=t_{0}}^{t=t_{f}} L[x,\dot{x},t]\,dt = 0.
\end{equation}
The lovely thing about an action principle, when we 
look at it this way, is that it provides a recipe
for constructing new equations of motion that are 
the result of replacing $L[x,\dot{x},t]$ by some
conveniently chosen approximation 
$L_{\rm app}[x,\dot{x},t]
\approx L[x,\dot{x},t]$.

If we make sure that our system of mathematics is
{\em large enough} to capably handle all useful
kinds of approximation, then --- once we have
formalized these in the abstract --- we will
obtain an entire new system of generalized
quantum dynamics.

So, the goal is to replace {\em exact} action principles
by {\em approximate} action principles and so obtain
entire {\em approximate dynamical systems}. Then we
look at these, study and classify them, the better
to understand their particular merits and deficiencies.

\section{Decorrelation as a standard approximation}
Practical approximations are designed to leave some
effect out to make things simple. In quantum theory
the one {\em generic effect} which makes the theory
hard to calculate with, and vizualize, is quantum
correlation and quantum entanglement. In the
theory of approximate dynamical systems we
use some simple tricks to suppress this
effect and simplify things.

Here is one simple semi--classical example
\begin{equation}
\langle\frac{\hat{p}^{2}}{2m} + k \hat{x}^{2}\rangle
\approx \frac{\langle\hat{p}\rangle^{2}}{2m}
+ k \langle \hat{x} \rangle^{2},
\end{equation}
where the quantum expectation is replaced by its
semi--classical counterpart. A familiar many--body
example is the Hartree approximation
\begin{eqnarray}
\lefteqn{\int \psi^{*}({\bf x}_{1},{\bf x}_{2})
       V(|{\bf x}_{1} - {\bf x}_{2}|)
         \psi({\bf x}_{1},{\bf x}_{2})
\,d^{3}{\bf x}_{1} d^{3}{\bf x}_{2} 
\approx } \nonumber \\
& & \hspace{3cm}
\int \psi^{*}_{1}({\bf x}_{1})\psi^{*}_{2}({\bf x}_{2})
       V(|{\bf x}_{1} - {\bf x}_{2}|)
 \psi_{1}({\bf x}_{1})\psi_{2}({\bf x}_{2})
\,d^{3}{\bf x}_{1} d^{3}{\bf x}_{2},
\end{eqnarray}
where $\psi({\bf x}_{1},{\bf x}_{2})$ has
been replaced by a factorized pair of
wave--functions.
 
In both cases we neglect correlations, or enforce
disentanglement, and so modify the degree of the
original expression in $\psi$ and $\psi^{*}$. It 
is this modification of degree which is the
cause of {\em induced nonlinearity}, as we
now see with a simple example. 

\section{Classical Schr\"{o}dinger equation}
The easiest way to express the correspondence between
classical and quantum physics is via the Ehrenfest
theorem. Starting with the quantum equations:
\begin{equation}
\label{exact}
\frac{d\langle\hat{p}\rangle}{dt} = 
-\langle H_{x}(\hat{x},\hat{p}) \rangle,
\;\;\mbox{ and }\;\;
\frac{d\langle\hat{x}\rangle}{dt} = 
+\langle H_{p}(\hat{x},\hat{p}) \rangle;
\end{equation}
we introduce the obvious semi--classical approximation:
\begin{equation}
\label{semiclass}
\langle H_{x}(\hat{x},\hat{p})\rangle \approx
H_{x}(\langle\hat{x}\rangle,\langle\hat{p}\rangle),
\;\;\mbox{ and }\;\;
\langle H_{p}(\hat{x},\hat{p})\rangle \approx
H_{p}(\langle\hat{x}\rangle,\langle\hat{p}\rangle);
\end{equation}
and thus obtain the {\em approximate\/} equations:
\begin{equation}
\label{approx0}
\frac{d\langle\hat{p}\rangle}{dt} \approx
- H_{q}(\langle\hat{x}\rangle,\langle\hat{p}\rangle),
\;\;\mbox{ and }\;\;
\frac{d\langle\hat{x}\rangle}{dt} \approx
+H_{p}(\langle\hat{x}\rangle,\langle\hat{p}\rangle).
\end{equation}
If we now take these as defining a new dynamical
system (i.e. we replace $\approx$ by $ = $) then 
our equations reduce to those of Hamilton,
\begin{equation}
\label{approx}
\frac{dP}{dt} =
- H_{x}(X,P),
\;\;\mbox{ and }\;\;
\frac{dX}{dt} =
+H_{p}(X,P);
\end{equation}
where we make the obvious identification:
\begin{equation}
\label{classcoord}
X(t) = \langle\hat{x}\rangle(t)
\;\;\mbox{ and }\;\;
  P(t) = \langle\hat{p}\rangle(t).
\end{equation}
In taking these steps one reduces the quantum
problem to a classical problem, in a manner
that ignores certain features of the full
quantum treatment.

Now let us apply this analysis of the Ehrenfest theorem,
as a decorrelation approximation, at the general level 
of the exact quantum action principle
\begin{equation}
\label{aexact}
\frac{\delta }{\delta \psi^{*}}
\int i\hbar\langle\psi|\frac{d}{dt}|\psi\rangle
- \langle\psi|\hat{H}(\hat{x},\hat{p})|\psi\rangle\,dt = 0.
\end{equation}
Taking variations with this we obtain
\begin{equation}
i\hbar\frac{d}{dt}|\psi\rangle
= \hat{H}(\hat{x},\hat{p})|\psi\rangle,
\end{equation}
as the general equation of motion. However, we could 
just as
well substitute 
\begin{equation}
\label{semiapp}
\langle\psi|\hat{H}(\hat{x},\hat{p})|\psi\rangle
\approx 
\langle\psi|H (\langle\hat{x}\rangle,
               \langle\hat{p}\rangle)|\psi\rangle.
\end{equation}
for the energy expectation, and so obtain directly
a decorrelated {\em classical wave--equation}.

There is, however, a minor subtlelty to carrying
out this program. In (\ref{semiapp}) it is not
guaranteed that the action principle remains
invariant to a re--normalization of $\psi$.
Obviously we want to retain that freedom 
to adjust and preserve normalization. To
overcome this difficulty we rescale all 
coordinate expectations as:
\begin{equation}
\label{homo0}
\langle\hat{x}\rangle =
\langle\psi|\hat{x}|\psi\rangle/n
\;\;\mbox{and}\;\;
  \langle\hat{p}\rangle =
\langle\psi|\hat{p}|\psi\rangle/n,
\end{equation}
where $n =\langle\psi|\psi\rangle$. Calculating
variational derivatives we find
\begin{equation}
\label{homop}
\frac{\delta\langle\hat{x}\rangle}{\delta\psi^{*}}
 = n^{-1}(\hat{x} - \langle\hat{x}\rangle)|\psi\rangle,
\;\;\mbox{and}\;\;
\frac{\delta\langle\hat{p}\rangle}{\delta\psi^{*}}
 = n^{-1}(\hat{p} - \langle\hat{p}\rangle)|\psi\rangle.
\end{equation}
Invoking now the approximate action principle 
\begin{equation}
\label{aapprox}
\delta \int i\hbar\langle\psi|\frac{d}{dt}|\psi\rangle
- \langle\psi|H
(\langle\hat{x}\rangle,\langle\hat{p}\rangle)
|\psi\rangle\,dt = 0,
\end{equation}
we use the chain rule
$$\frac{\delta }
{\delta\psi^{*}}\left[
\langle\psi|
H(\langle\hat{x}\rangle,\langle\hat{p}\rangle)|\psi\rangle\right]
= 
H(\langle\hat{x}\rangle,\langle\hat{p}\rangle)
\frac{\delta n}
 {\delta\psi^{*}} +
 H_{x}(\langle\hat{x}\rangle,\langle\hat{p}\rangle)
\frac{\delta \langle\hat{x}\rangle}
{\delta\psi^{*}} +
  H_{p}(\langle\hat{x}\rangle,\langle\hat{p}\rangle)
\frac{\delta \langle\hat{p}\rangle}
{\delta\psi^{*}},$$
to obtain the approximate equation of motion
\begin{equation}
i\hbar\frac{d}{dt}|\psi\rangle
=\left\{
H(\langle\hat{x}\rangle,\langle\hat{p}\rangle)\hat{1}
+
H_{x}(\langle\hat{x}\rangle,\langle\hat{p}\rangle)
(\hat{x} - \langle\hat{x}\rangle)
+
H_{p}(\langle\hat{x}\rangle,\langle\hat{p}\rangle)
(\hat{p} - \langle\hat{p}\rangle)\right\}
|\psi\rangle.
\end{equation}
This is the {\em classical Schr\"{o}dinger equation},
which recovers the Ehrenfest equations of motion in
classical form. It propagates wave--packets
neglecting {\em dispersion} and {\em correlation}.
The result is that they bounce off barriers and 
the like just like classical particles.

One can construct exact solutions of the above
nonlinear integrodifferential equation. To do
this we first solve the classical problem to
find $X(t)$ and $P(t)$. Next we take any
wavefunction $\psi_{0}(x)$ having both
position and momentum expectation
values equal to zero. Then we
form the time--dependent wavefunction
\begin{equation}
\label{solutions}
\psi(x,t) = 
e^{\frac{i}{\hbar}\int_{t_{0}}^{t} L\,d\tau}
e^{-iP(t)X(t)/2\hbar}
e^{iP(t)x/\hbar}\psi_{0}(x - X(t)),
\end{equation}
where the exact classical action
\begin{equation}
\label{phase}
\int_{t_{0}}^{t} L\,d\tau = 
\int_{t_{0}}^{t}
\left(\frac{P\dot{X} - X\dot{P}}{2}\right)
- H(X,P) \,d\tau 
\end{equation}
appears as the leading phase factor (showing
that the Feynmann--Dirac correspondence
is semi--classically exact). This argument 
can be made constructive, but it is
much easier to verify by substitution. Alternatively,
given a theory of nonlinear propagators one can set 
up these equations on a computer and solve them 
directly to verify this general solution.

\section{Physical interpretation}
When dealing with approximate dynamical
systems we must remember that linearity is vital 
to the Copenhagen interpretation. However, we 
use nonlinear wave--equations all the time 
in physics. To interpret them we adopt a
computational algorithm viewpoint.

We have an exact theory, and quantities that we wish
to calculate --- e.g. eigenvalues, stationary and 
time--dependent wavefunctions, expectation values, 
transition probabilities etc. These we could 
calculate {\em exactly} or {\em approximately}.

Either way we {\em can} apply a physical interpretation 
that presupposes linearity as an exact property of
nature. To the approximately computed, i.e.
nonlinearly evolved, physical quantities
we apply the Copenhagen interpretation
--- on the understanding that there
is supposed to be an error in our
treatment somewhere.\footnote{Of course, ultimately 
the matter of what is correct rests with experiment.}

For instance, in solving our approximate
classical equations we have no need of
$\hbar$, nor any specific wavefunction.
The solution of the reduced, and thus
simplified, problem requires only the
initial expectation values. It is thus 
an approximate method for computing 
$\langle\hat{x}\rangle(t)$ and  
$\langle\hat{p}\rangle(t)$. 
The errors committed are identical, both 
numerically and conceptually, to those
of the familiar Hamiltonian dynamics. 
Even so, it is pretty useful.
One could say the existence of
the classical Schr\"{o}dinger
equation, as an excellent
semiclassical approximation, 
{\em explains} why classical 
dynamics fooled us for 300 years!

\end{document}